\documentclass[12pt]{article}
\usepackage{amsmath,amssymb,latexsym}
\usepackage{times}
\usepackage{graphicx}

\newtheorem{theorem}{Theorem}[section]
\newtheorem{remark}{Remark}[section]

\numberwithin{equation}{section}

\begin{document}
\title{Letter Representations of $m\times n\times p$ Proper Arrays}
\date{\today}
\author{Jocelyn Quaintance\\
quaintan@temple.edu}
\maketitle
\begin{abstract}
\noindent Let $m\neq n$.  An $m\times n\times p$ {\it proper array} is a three-dimensional rectangular array composed of directed cubes that obeys certain constraints.  Because
of these constraints, the $m\times n\times p$ proper arrays may be classified via a schema in which each $m\times n\times p$ proper array is associated with a particular $m\times n$ planar face.  By
representing each connected component present in the $m\times n$ planar face with a distinct letter, an $m\times n$ array of letters is formed.  This $m\times n$ array of letters is the {\it letter
representation} associated with the
$m\times n\times p$ proper array.  The main result of this paper involves the enumeration of all $m\times n$ letter representations modulo symmetry, where the symmetry is derived from the group $D_2 =
C_2\times C_2$ acting on the set of letter representations.   The enumeration is achieved by forming a linear combination of four exponential generating functions, each of which is derived from a
particular symmetry operation.  This linear combination counts the number of partitions of the set of $m\times n$ letter representations that are inequivalent under $D_2$.
\end{abstract}
\newpage
\section*{Introduction}
\subsection {History of the Problem: Two-Dimensional Proper Arrays}
\noindent This article is an extension of the author's work involving the enumeration of two-dimensional proper arrays [1],[2],[3].  In both the two-dimensional and three-dimensional situations, the
basic building block of a proper array is a {\bf directed cube}, where a directed cube is a cube with five faces that containing indentations while the sixth face has a cylindrical plug or {\bf
connector}.
\begin{center}
\centerline{\includegraphics[width=3.0cm, height=3.0cm]{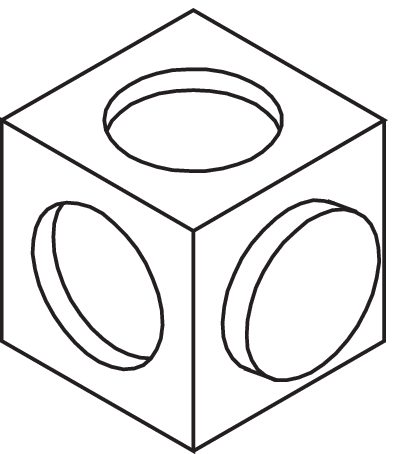}}
\end{center}
\begin{center}
\small{ Figure 1: A directed cube.  In this illustration, the connector extends from the right face while the other five faces contain circular indentations.}
\end{center}
\noindent   Rectangular arrays can be built by fitting the connector of one directed cube into an indentation of another.  If the ensuing rectangular array, which consists $m$ rows and $n$ columnns, obeys
certain connectivity constraints, it is said to be an {\boldmath $m\times n$} {\bf proper array} [1],[2],[3].  The $m\times n$ proper array may be classified by its rightmost column of directed cubes.  The
right edge structure contains all of the path-connected subsets that appear in the $m\times n$ proper array and all of the connectors that extend from this edge.\\\\ 
\noindent The right edge structure does not uniquely determine the number $m\times n$ proper arrays since two different $m\times n$ proper arrays may have the same right edge structure. 
However, the right edge structure is the state for a transition matrix, {\boldmath $M_m$}, which for a fixed $m$ and arbitrary $n$, describes how an $m\times n$
proper array evolves into an $m\times (n+1)$ proper array [1], [2], [12].  Note that {\boldmath $M_m$} is an $r\times r$ matrix, where $r$ is the number, modulo symmetry, of possible
rightmost columns associated with $m\times n$ proper arrays [1],[2],[12].  Let {\boldmath $S_n$} be the $1\times r$ vector $(s_{n,i})$, where $s_{n,i}$ counts the number of $m\times n$ proper
arrays, modulo symmetry, having the $i$th kind of right edge.  We find {\boldmath $S_nM_m$} $=$ {\boldmath $S_{n+1}$}.  For example, in the case of $2\times n$ proper arrays, there are four possible right
edge structures.  Hence, {\boldmath $M_2$} is a $4\times 4$ matrix.  Let {\boldmath $S_1$} $ = [1,\,\,2, \,\, 0,\,\, 1]$, where the entries of {\boldmath $S_1$} count the number of $2\times 1$ proper
arrays, with a particular right edge structure.  Note, there are $1 + 2 + 0 + 1 = 4$ distinct $2\times 1$ proper arrays.  Then {\boldmath $S_1M_2$} $=$ {\boldmath $S_2$} $ =
[2,\,\,8,\,\,0,\,\, 3]$, where the entries of {\boldmath
$S_2$} count the number of $2\times 2$ proper arrays with a particular right edge structure.  Note, there are
$2 + 8 + 0 + 3 = 13$ distinct $2\times 2$ proper arrays.  Thus, the right edge structure of an $m\times n$ proper array can determine the total number, modulo symmetry, of
$m\times n$ proper arrays [1],[2],[12].\\\\ 
\noindent For a fixed $n$, the author wrote a computer program which constructs {\boldmath $M_m$}.  In the process of writing this program, the author realized that one of the main research questions
associated with the enumeration of $m\times n$ proper is determining, for any given $n$, the size of $r$, where $r$ is the number of right edge structures, modulo symmetry, associated with $m\times n$
proper arrays.  By having such a formula, the author and various computer scientists [5] would be able to determine the amount of computer time and memory necessary for the construction of
{\boldmath $M_m$}.  This formula consists of a linear combination of four ordinary generating functions, each of which enumerate a particular subset of right edge configurations [1], [3].  Thus, for the
two-dimensional situation, the basis size for {\boldmath $M_m$} is completely determined.
\begin{center}
\begin{tabular}{||l|r|r|r|r|r|r||}
\hline
m & 1 & 2 & 3 & 4 & 5 & 6\\
\hline
{\boldmath $M_m$} & 1 & 4 & 16 & 72 & 364 & 1916\\
\hline
\end{tabular}
\end{center}
\begin{center}
\small{Table 1:  The basis size for the cellular automata {\boldmath $M_m$}.}
\end{center}

\subsection{Three-Dimensional Proper Arrays}
\noindent Next, the author decided to enumerate three-dimensional proper arrays.  An {\boldmath $m\times n\times p$} {\bf proper array} is an $m\times n\times p$ array of directed cubes that obey certain
constraints.  These constraints allow the set of $m\times n\times p$ proper arrays to be classified via a schema in which each $m\times n\times p$ proper array is associated with
a particular $m\times n$ planar face, called the {\bf preferred face}.  The preferred face of the $m\times n\times p$ proper array contains all the path-connected subsets that appear in the
$m\times n\times p$ proper array and the number of connectors that extend from the preferred face. \\\\ 
\noindent After defining the notion of an $m\times n\times p$ proper array, it was an easy exercise to adapt the computer program
which constructed {\boldmath $M_m$} into one which constructed {\boldmath $M_{m\times n}$}, where {\boldmath $M_{m\times n}$} is the transition matrix that describes how an $m\times n\times p$ proper array
evolves into an $m\times n\times (p+1)$ proper array.  Note that {\boldmath $M_{m\times n}$} is an $r\times r$ matrix, where $r$ counts, modulo symmetry, the number of $m\times n$ preferred faces
associated with $m\times n\times p$ proper arrays.  The most important adaptation involved the varying numbers of symmetry maps.  Unlike
$m\times n$ proper arrays, which have only one symmetry map, $m\times n\times p$ proper arrays have four or eight symmetry maps.  More precisely, if
$m\neq n$, an $m\times n\times p$ proper array possesses the rectangular symmetry of $D_2$.  However, when $m = n$, an $m\times m\times p$ proper array has $D_4$ symmetry.  Each of these symmetry maps
must take the preferred face onto itself.\\\\
\noindent Once again, we would like to find a formula that predicts $r$, the basis size of {\boldmath $M_{m\times n}$}.  The main result of this paper is a first step to constructing such a forumla, when
$m\neq n$.  In particular, Theorem 2.1 provides a formula, in the form of a linear combination of four exponential generating functions, for calculating, modulo $D_2$ symmetry, the number
of preferred face structures associated with $m\times n\times p$ proper arrays.  The techniques used to derive the generating functions are similar in nature to the techniques used by
Yoshinaga and Mori [7] and David Branson [8].  The second result provides generating functions which count the number of $m\times n\times p$ proper arrays whose associated
$m\times n$ preferred face is fixed by the four symmetry maps of $D_2$.
\subsection{Connections to Percolation Theory}
\noindent Before continuing, we would like to provide a physical interpretation for an $m\times n\times p$ proper array.  If the center of each directed cube is represented as a point of
$Z^3$, with the connectors providing a network of open paths between the various vertices, then the $m\times n\times p$ proper array is a model of bond percolation [10 Section 1], [11
P.16].  Furthermore, this bond percolation has the following property, namely, there exists an open cluster that allows water to flow from the back face of the $m\times n\times p$ array to its front face
[10 P.2].  Moreover, since there is a relationship between bond percolation and the Ising Model [10 P.8], an open question involves the exact nature of the connection between $m\times n\times p$ proper
arrays and Ising Models.  At this point in time, we have not investigated this connection between the Ising Model and the three-dimensional proper arrays and leave it as an open question for future
research. 
\section{Proper Arrays and Letter Representations}
\noindent We are now ready to give the definition of an $m\times n\times p$ proper array.\\\\
\noindent {\bf Definition:} A (connected) {\bf component} of an array of directed cubes is a (path) connected subset of the array.\\\\
\noindent {\bf Definition:} Let $A$ represent an $m\times n\times p$ array of directed cubes.  Orient $A$ so that one planar $m\times n$ face has a center at $(0,1,0)$ and is perpendicular to the $y$ axis. 
This $m\times n$ planar face is the {\bf preferred face} of $A$.\\\\
\newpage
\noindent {\bf Definition:} An {\boldmath $m\times n\times p$} {\bf proper array} is a three-dimensional array of directed cubes that obeys the following two conditions:
\begin{itemize}
\item[1.]  The array is {\bf 5-way flat}.  That is, only the preferred face has outward pointing connectors.
\begin{center}
\includegraphics[width=1.75in, height=2.0in]{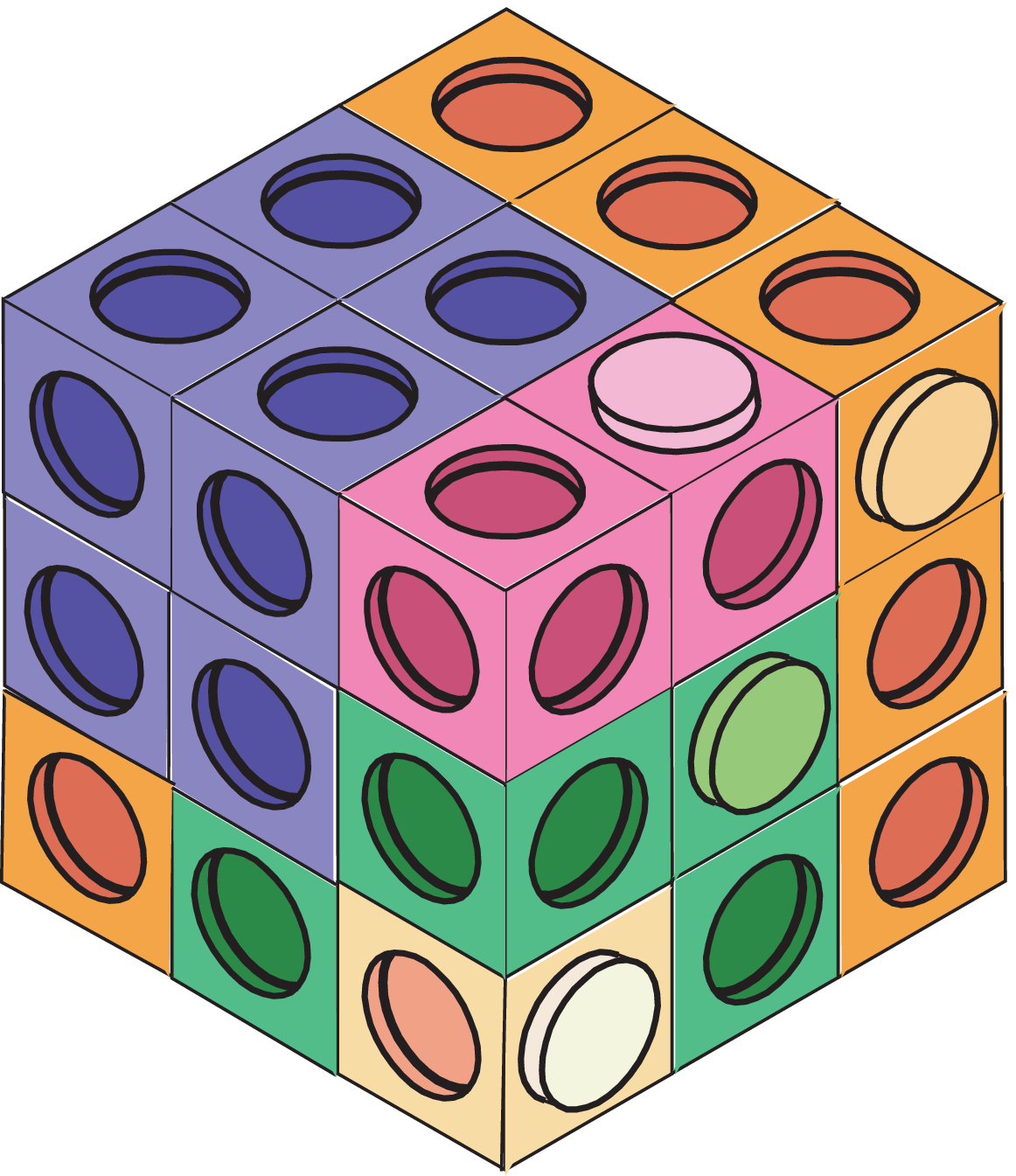}
\end{center}
\begin{center}  
\small{Figure 1.1: This $3\times 3\times 3$ array of directed cubes is not 5-way flat since the pink component has a conncetor that extends from the top face.}
\end{center}
\item[2.]  We define an {\bf island} to be a component that does not reach the preferred face. The array \underline{can not} have islands.
\begin{center}
\includegraphics[width=1.75in, height=2.0in]{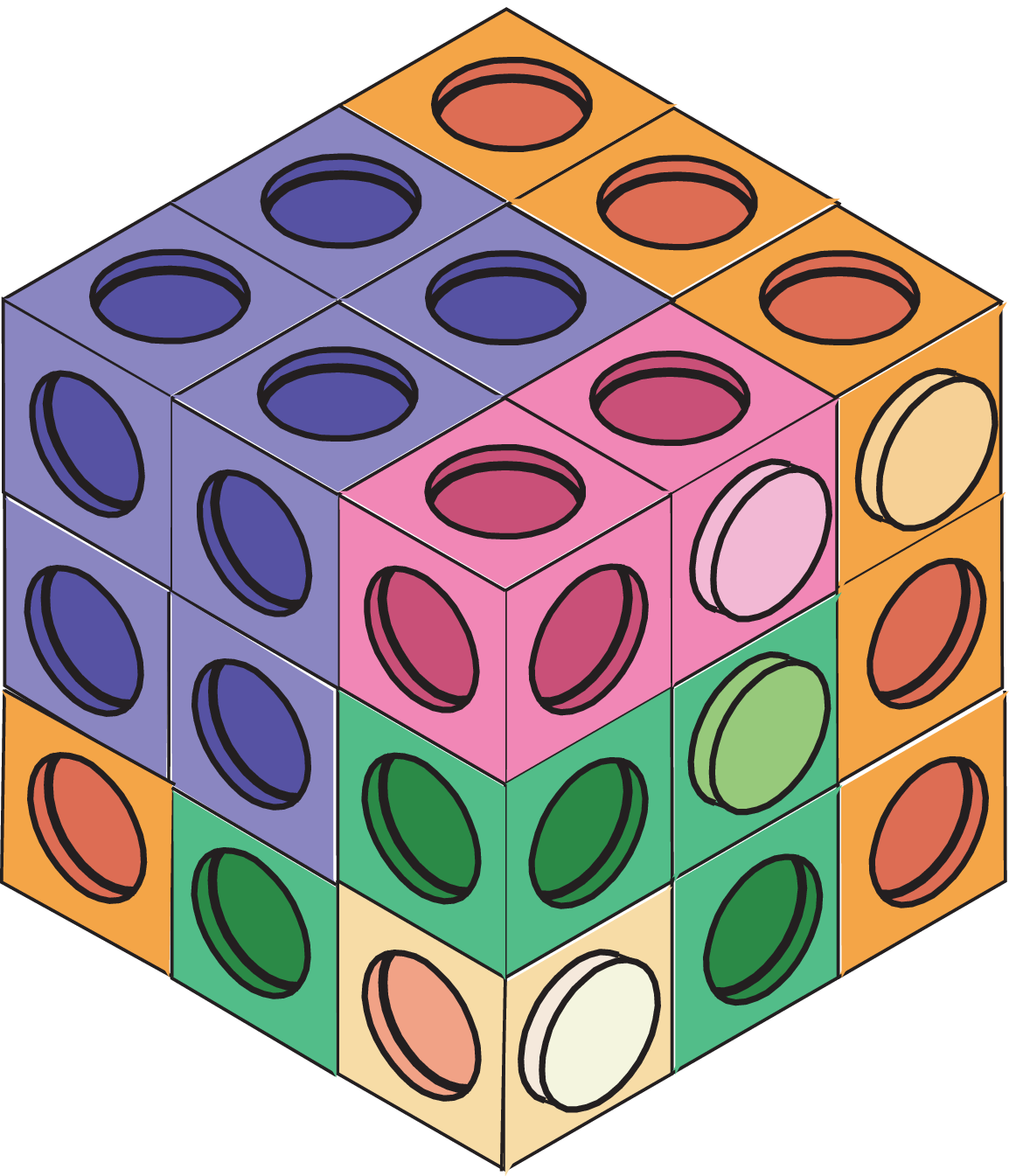}
\end{center}
\begin{center} 
\small{Figure 1.2: The purple component forms an island since it does not reach the preferred face.}
\end{center}
\end{itemize}
\noindent An $m\times n\times p$ proper array may be associated with its preferred face, since the preferred face corresponds to the state of a transition matrix.  In
particular, by letting each distinct component be denoted by a distinct letter, and any outward pointing connectors be depicted as circles around the letters, the geometry of the preferred face
is recorded as $m\times n$ array of circled letters.  This $m\times n$ array of circled letters is called the {\bf word
representation} of the $m\times n\times p$ proper array.  It encodes the numbers of connected components that appear in the $m\times n\times p$ proper array and the connector that extend from the $m\times
n\times p$ proper array.  By deleting the circles in the word representation, the {\bf letter representation} of the $m\times n\times p$ proper array is formed.
\begin{center}
\includegraphics[width=1.75in, height=2.5in]{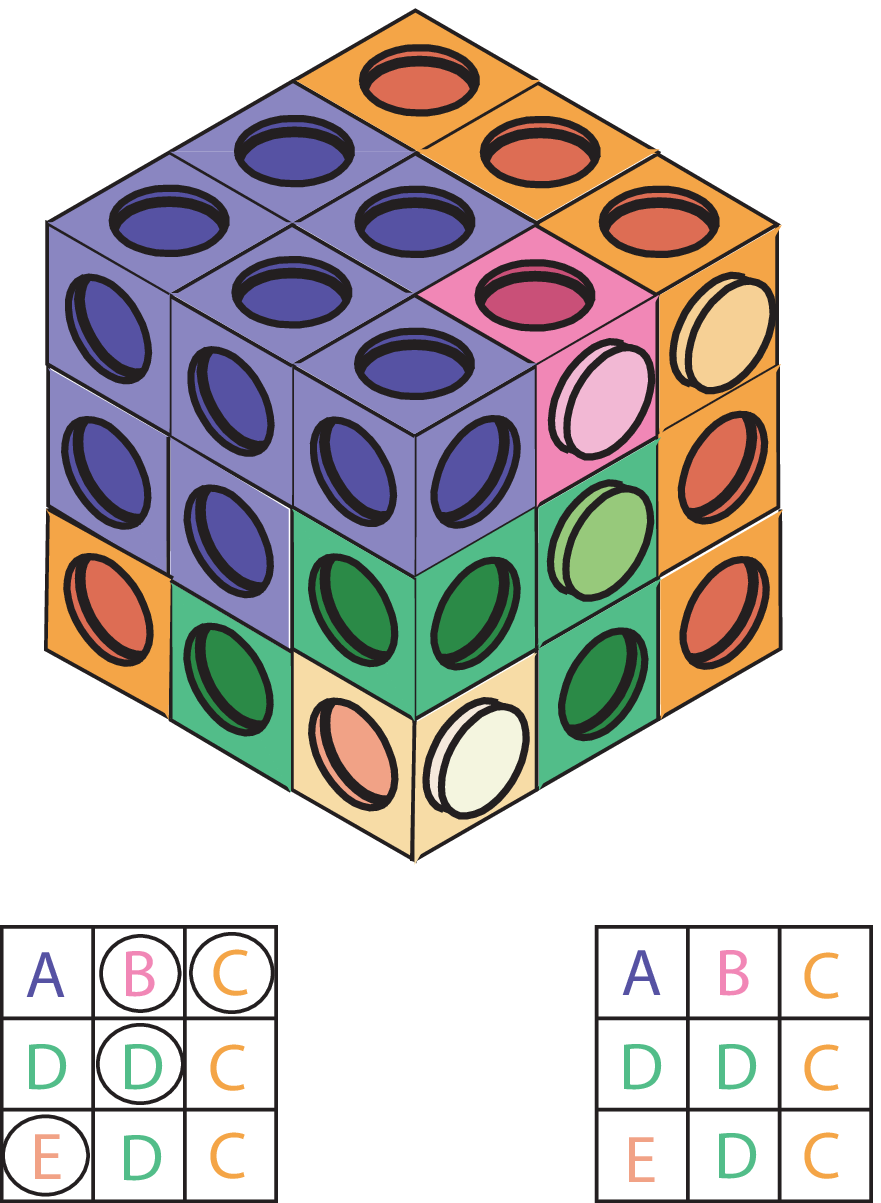}
\end{center}
\begin{center} 
\small{Figure 1.3: An example of a $3\times 3\times 3$ proper array.  At the bottom left corner, the word representation is given.  This word
representation records the components as colored letters and outward pointing connectors as circles present.  At the bottom right corner, the letter
representation is given.}
\end{center}  
\begin{remark}
Letter Representations are partitions of the preferred face.  Hence, two letter representations are regarded as the same if they have the same components but their letters are different.  In other words,
the letters are just labels for the components and the labels are unimportant.  For example, take the letter representation provided in Figure 1.3.  If we replace the $A$ with an $X$, the resulting
letter representation is the same as the orginal letter representation.
\end{remark} 
\noindent The definition of an $m\times n\times p$ proper array is motivated by the desire to ensure that a particular $m\times n$ face of the proper array corresponds to the state of the  transition matrix
[1],[2],[12].  In order to see how the definition of proper array guarantees this correspondence, we will take a moment to describe how the transition matrix constructs proper arrays.  Fix $m$ and $n$. 
The matrix multiplication {\boldmath $S_{p-1}M_{m\times n}$} corresponds to attaching an
$m\times n$ layer of directed cubes onto a previous constructed $m\times n\times (p-1)$ proper array.  The concept of 5-way flat ensures that this conncection could only occur on the preferred face. 
The concept of no islands ensures there are no connected subsets of the $m\times n\times (p-1)$ proper array that do not reach its preferred face.  Hence, all the connected subsets of the $m\times n\times
(p-1)$ proper array will be affected by the attachment of the new $m\times n$ layer and will percolate throughout the newly constructed $m\times n\times p$ proper array.  Since this percolation occurs, we
are able to record the state of construction by simply looking at the geometry of the $m\times n$ preferred face of the newly build $m\times n\times p$ proper array.  Without the conditions of 5-way flat
and no islands, this correspondence between preferred face and construction state would not exist.

\section{Enumerating Letter Representations}
\noindent  For the remainder of this paper, we can ignore the three-dimensional context provided by the proper array and work in the two-dimensional
setting of the letter representation.\\\\
\noindent Our goal is to fix $m$ and $n$, and enumerate, modulo symmetry, all letter representations
associated with $m\times n\times p$ proper arrays.  If $m\neq n$, the symmetry equivalence is determined by $D_2$.  
\begin{center}
\includegraphics[width=3.0in, height=3.5in]{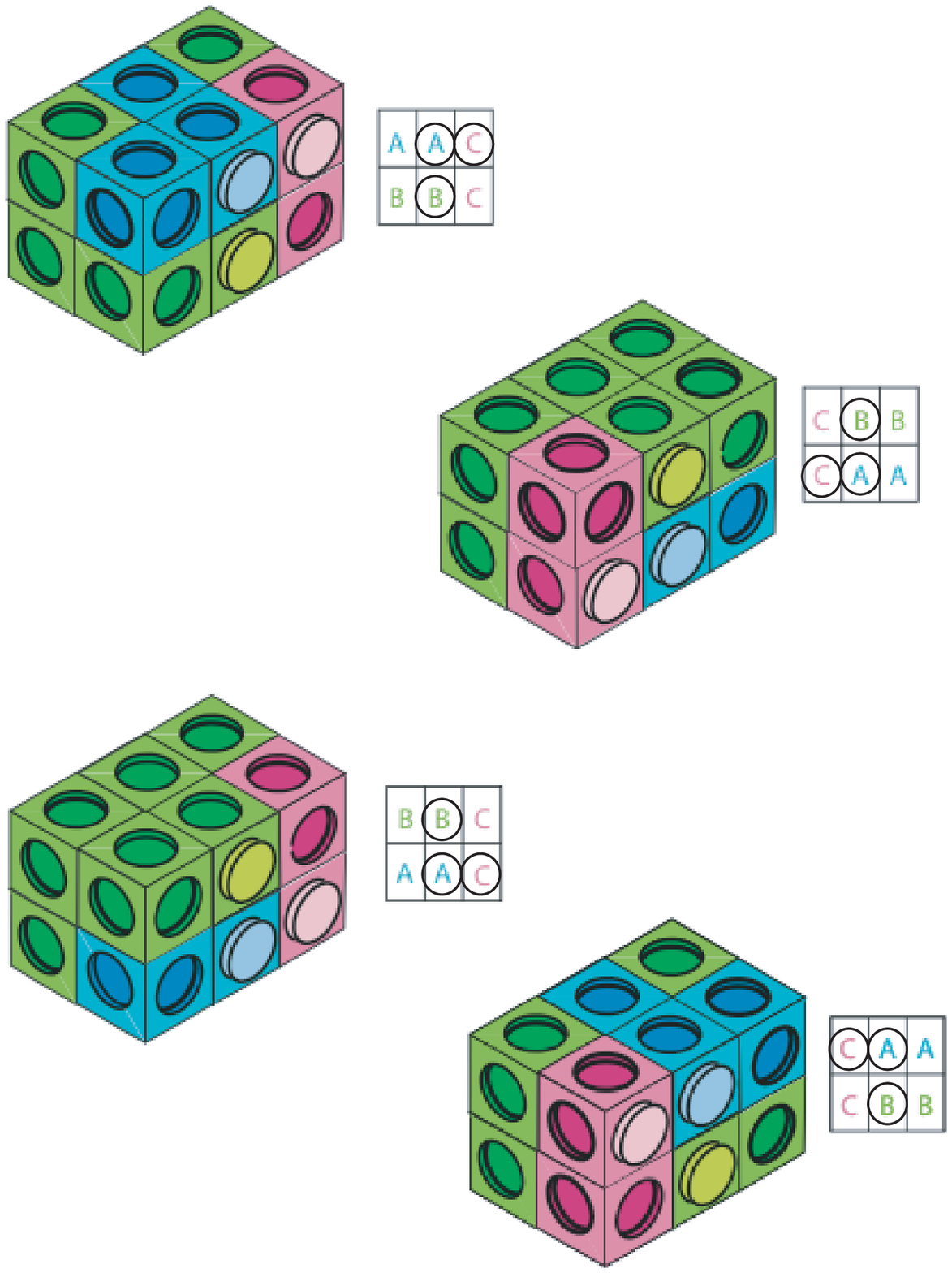}
\end{center}
\begin{center}
\small{Figure 2.1: The four equivalent versions of a $2\times 3\times 2$ array.  Let the diagram in the upper left corner represent the orginal array.  Then by reading in a clockwise manner, the next
diagram depicts a $180^{\circ}$ rotation of the orginal array.  The bottom right corner depicts a vertical reflection of the orginal array.  The remaining diagram depicts a horizontal reflection of the
orginal array.  These four symmetry images are only counted once in the enumeration procedure.} 
\end{center}

\noindent {\bf Definition:}  A {\bf {\boldmath $m\times n$} letter array} is an $m\times n$ arrangement of letters which obeys the following condition:  Each distinct letter labels an element in the
partition of the $m\times n$ array of squares.  The comments of Remark 1.1 applly to the set of $m\times n$ letter arrays.\\\\
\noindent  It is an easy exercise to show that $B_{m,n}$, the cardinality of the entire set of
$m\times n$ letter arrays, equals to the $mn$th Bell Number.  Define $L_{m,n}$ to be the number of
$m\times n$ letter representations modulo $D_2$ symmetry.  Note that $L_{m,n}$ is a lower bound for the basis of the transition matrix {\boldmath $M_{m\times n}$}. We
calculate $L_{m,n}$ as follows.  
\begin{itemize}
\item[1.] Let $S_{m,n}$ count the $m\times n$ letter arrays that are fixed via horizontal reflection, vertical reflection, and $180^{\circ}$ degree rotation.  
\item[2.] Let $H_{m,n}$ count the $m\times n$ letter arrays fixed via horizontal reflection.  Then, $H_{m,n} - S_{m,n}$ counts the $m\times n$ letter arrays that are fixed only by horizontal
reflection.  
\item[3.] Let $V_{m,n}$ count the $m\times n$ letter arrays fixed via vertical reflection.  Then, $V_{m,n} - S_{m,n}$ counts the $m\times n$ letter arrays that are fixed only by vertical
reflection. 
\item[4.] Let $R_{m,n}$ count the $m\times n$ letter arrays fixed via $180^{\circ}$ rotation.  Then, $R_{m,n} - S_{m,n}$ counts the $m\times n$ letter arrays that are fixed only by
rotation.  
\item[5.] Let $C_{m,n} = B_{m,n}-(H_{m,n}-S_{m,n})-(V_{m,n}-S_{m,n})-(R_{m,n}-S_{m,n})-S_{m,n}$.  Then, $C_{m,n}$ counts the $m\times n$ letter arrays that are not fixed by any symmetry
transformation.
\end{itemize}
\begin{theorem}
Let $L_{m,n}, C_{m,n}, B_{m,n}, V_{m,n}, H_{m,n}, R_{m,n},$ and $S_{m,n}$ be as previously defined.  Then 
\[ L_{m,n} =  \frac{B_{m,n}+H_{m,n}+V_{m,n}+R_{m,n}}{4}\]
\end{theorem}
{\bf Proof of Theorem 2.1:} To calculate the letter representations associated with $m\times n\times p$ proper arrays modulo symmetry, we first determine whether a given $m\times n$ letter array, called
$A$, is fixed via any of the four symmetry transformations.  If $A$ is not fixed by any symmetry, it has four equivalent images.  However, if $A$ is fixed
under a symmetry transformation, it has at most two symmetry images.  It follows that  
\begin{align*}
L_{m,n} &= \frac{C_{m,n}}{4} + \frac{(H_{m,n}-S_{m,n})+(R_{m,n}-S_{m,n})+(V_{m,n}-S_{m,n})}{2} + S_{m,n}\\
&= \frac{B_{m,n}+H_{m,n}+V_{m,n}+R_{m,n}}{4} \qquad \Box
\end{align*}
\begin{remark}
Theorem 2.1 can be considered to be an immediate consequence of Burnside's Lemma.
\end{remark}
\subsection{Numerical Example}
\noindent In order to understand how the information provided by Theorem 2.1 helps provide a lower bound on the basis size of the transition matrix, look at the following example.  Let $m = 3$ and $n =
1$.  The goal is to compute the size of the transition matrix associated with the $3\times 1\times p$ proper arrays.  It can be shown that the transition matrix is a $16\times 16$ matrix [1],[2].  Thus,
there are 16 possible $3\times 1$ {\it word} representations modulo symmetry.  Theorem 2.1 simply counts the {\it letter}
representations modulo symmetry.  Ignoring symmetry, there are $B_{3,1} =  5$ letter representations.  By direct calculation, we find that $H_{3,1} = 3,  V_{3,1} = 5,$ and  $R_{3,1} = 3$.  Theorem 2.1
implies that the number of $3\times 1$ letter representations modulo symmetry is $\frac{5 + 3 + 5 + 3}{4} = 4$.  Thus, the transition matix associated with the $3\times 1\times p$ proper arrays is at
least a $4\times 4$ matrix.  The goal of our research is to obtain a formula that calculates the actual basis size of the transition matrix.  As this example demonstrates, the theorems in
this paper provide not actual basis size, but a lower bound.  
\section{Generating Function for Horizontal/Vertical Symmetry}
\noindent In order to use Theorem 2.1, we need to find generating functions for $H_{m,n}$, $V_{m,n}$, and $R_{m,n}$, These generating functions are obtained by dividing the
$m\times n$ array into small sections determined by the symmetry transformation.  For instance, in the case of horizontal reflection, the array is divided into two halves.  If
two or more transformations are applied, the array is divided into four quadrants.  In either case, we can arbitrarily fill one of the halves/quarters with any arrangement of letters and then use
symmetry to fill the remaining half/quarters of the $m\times n$ array.  The key in this technique is to carefully divide the array around the row and/or column that may be fixed under a symmetry
transformation.  Hence, the generating functions depend on the parity of $m$ and $n$.\\\\
\noindent  We begin with the generating functions associated with $H_{m,n}$.  Due to symmetry, they are also the generating functions associated with $V_{m,n}$.  Before going further, we should recall the
following representations of the Bell Numbers and the Second Stirling Numbers [4],[7],[8].  
\begin{remark}
\[\sum_{n=0}^{\infty}B(n)\frac{x^n}{n!} = e^{e^x-1}\]
\end{remark}
\begin{remark}
\[\sum_{r=0}^{\infty}\left\{\begin{array}{c}r\\t\end{array}\right\}\frac{y^r}{r!} = \frac{(e^y-1)^t}{t!}\]
\end{remark}
\noindent Using the aforementioned subdivision technique, we find that
\[H_{2m,n} = \sum_{j=1}^{mn}\sum_{s=0}^{[\frac{j}{2}]}\frac{\left\{\begin{array}{c} mn\\j\end{array}\right\}2^{j-3s}j!}{s!(j-2s)!}\]
\noindent where
\begin{itemize}
\item[1.] $j$ counts the letters in the first $m$ rows.  We will refer to the first $m$ rows as the {\bf top half} of the array and the remaining $m$ rows as the {\bf bottom half} of
the array.
\item[2.] $s$ counts the number of these $j$ letters that are interchanged when reflected to the bottom half of the array.
\end{itemize}
\begin{center}
\hspace{3.5cm}
\includegraphics[width=6.0cm, height=4.2cm]{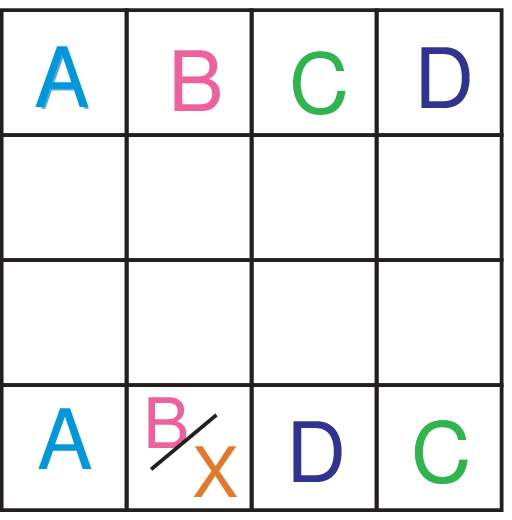}
\end{center}
\begin{center}
\vspace{-1.0cm}
\small{Figure 3.1: This figure demonstrates the three possiblities that occur to the letters that appear in the top half of the array.   A letter can reflect to itself; this is demonstrated by
$A$.  A letter can reflect to an entirely new letter; this is demonstrated by $B$ reflecting to $X$.  Otherwise, two letters, such as $C$ and $D$, can interchange position under reflection.} 
\end{center}
\begin{theorem}
The coefficient of $t^{mn}$ in the expansion of $\exp(\frac{1}{2}(e^t + 3)(e^t - 1))$ is $H_{2m,n}$.
\end{theorem}
\noindent {\bf Proof of Theorem 3.1:}
Let $k = mn$, then
\[H(t) = \sum_{k,j,s=0}^{\infty}
\frac{\left\{\begin{array}{c} k\\j \end{array}\right\}t^k2^{j-3s}j!}{s!(j-2s)!k!}\]
\noindent To obtain the desired result, use Remark 3.2 and sum over $s$ and $j$.  See the proof of Theorem 3.2 for a more elaborate example. \qquad $\Box$\\\\
\noindent Next we determine $H_{2m+1,n}$.  Finding $H_{2m+1,n}$ is more complicated since we
must contend with a middle row that is fixed by the reflection map. 
\[H_{2m+1,n} =n!\sum_{i=1}^{mn}\sum_{j=1}^{min(n,i)}\sum_{k=j}^{n}\sum_{s=0}^{[\frac{i-j}{s}]}\frac{\left\{\begin{array}{c} mn\\i\end{array}\right\}
\left\{\begin{array}{c} k\\j\end{array}\right\}B_{n-k}i!2^{i-j-3s}}{k!(n-k)!s!(i-j-2s)!}\]
\normalsize{\noindent where
\begin{enumerate}
\setlength{\itemsep}{-\itemsep}
\item[1.] $i$ counts the letters in the top half of the array.
\item[2.] $j$ counts the number of those $i$ letters that appear in the middle row. 
\item[3.] $k$ counts the squares in the middle row occupied by the $j$ letters.
\item[4.] $s$ counts the number of the $i-j$ letters that are interchanged when reflected to the bottom half of the array.
\item[5.] $B_{n-k}$ counts the $1\times (n-k)$ letter arrays.  These $1\times (n-k)$ letter arrays only appear in the middle row of array. 
\end{enumerate}
\begin{center}
\includegraphics[width=4.0cm, height=3.0cm]{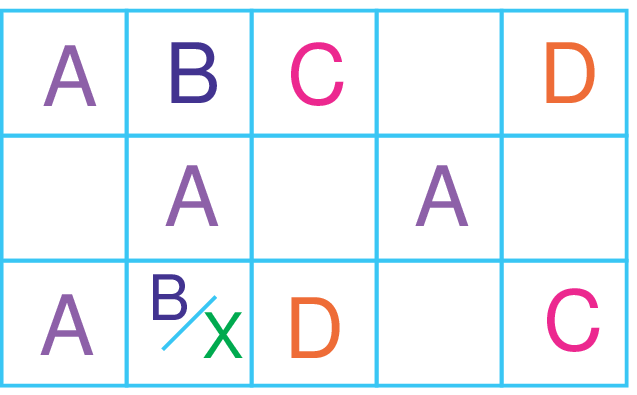}
\end{center}
\begin{center}
\small{Figure 3.2:  In this case, the middle row is fixed under horizontal reflection.  If a letter from the top half of the array occurs in the middle row, it is \underline{fixed}
under reflection.  Otherwise, it obeys the three possiblities discussed in the construction of $H_{2m,n}$.}
\end{center}
\begin{theorem}
The coefficient of $y^{mn}x^n$ in the expansion of $\exp(e^{y+x} + \frac{1}{2}e^{2y} - \frac{3}{2})$ is $H_{2m+1,n}$.
\end{theorem} 
\noindent {\bf Proof of Theorem 3.2}
Let $q= mn$ and define
\[H(x,y) = \sum_{\scriptstyle i,j,k,s,n,q = 0}^{\infty}\frac{\left\{\begin{array}{c} q\\i\end{array}\right\}\left\{\begin{array}{c}
k\\j\end{array}\right\}B_{n-k}i!2^{i-j-3s}x^ny^q}{k!(n-k)!s!(i-j-2s)!q!}\]
\noindent Using Remarks 3.1 and 3.2, we find 

\begin{align*}
H(x,y) &= \sum_{\scriptstyle i,j,s = 0}^{\infty} \frac{(e^y - 1)^i (e^x - 1)^j e^{e^x - 1} 2^{i - j - 3s}}{j!s!(i - j - 2s)!}\\
&= \sum_{\scriptstyle i,s =0}^{\infty}\frac{(e^y-1)^i(e^x+1)^{i-2s}e^{e^x-1}2^{-s}}{s!(i-2s)!}\\
&= \sum_{\scriptstyle s=0}^{\infty}\frac{(e^y-1)^{2s}e^{e^x-1 +(e^y-1)(e^x+1)}2^{-s}}{s!}\\
&= \exp(e^{y+x}+\frac{1}{2}e^{2y}-\frac{3}{2}) \qquad \Box
\end{align*}
\section{Generating Function for Rotational Symmetry}
\noindent Our next step is to find a generating function for the rotational symmetry.  Observe that\\ $R_{2m,2n} = H_{2m,2n}$, $R_{2m+1,2n} = H_{2n,2m+1},$ and $R_{2m,2n+1} = H_{2m,2n+1}$.  It remains to
compute $R_{2m+1, 2n+1}$.  In this situation, the $180^{\circ}$ rotation
fixes the central square of the array.  Our geometric subdivision must avoid this central square and takes the form of an
upper $L$ and a lower $L$.
\begin{center}
\includegraphics[width=3.5cm, height=3.5cm]{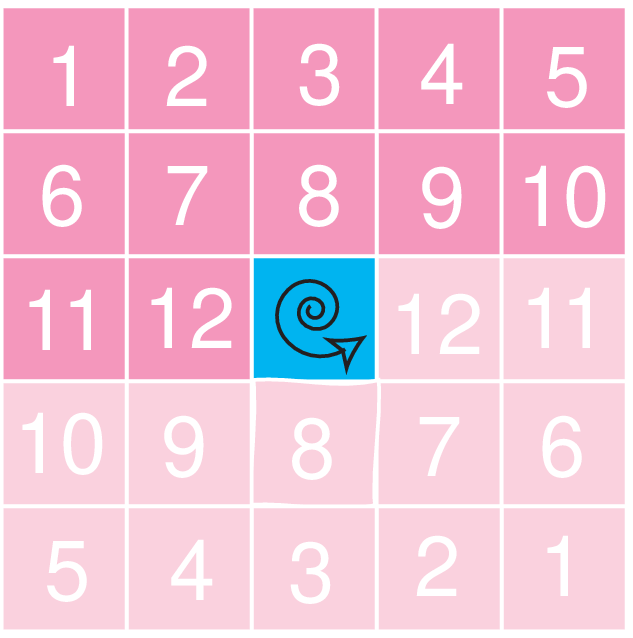}
\end{center}
\begin{center}
\small{Figure 4.1: Under $180^{\circ}$ rotation, the central square is fixed and the array is subdivided into two $L$ shapes.}
\end{center}
\noindent We find that 
\[R_{m,n} = \sum_{i=1}^{q}\sum_{s=0}^{[\frac{i}{2}]}\sum_{r=0}^{i-2s}\frac{\left\{\begin{array}{c} q\\i\end{array}\right\}i!(i-2s+1-r)}{s!2^sr!(i-2s-r)!}\]
\noindent where
\begin{enumerate}
\setlength{\itemsep}{-\itemsep}
\item[1.] $q = n[\frac{m}{2}] + [\frac{n}{2}]$.  
\item[2.] $i$ counts the letters in the top $L$ 
\item[3.] $s$ counts the number of the $i$ letters that are interchanged in the bottom $L$
\item[4.] $r$ counts the new letters that appear in the bottom $L$.  The central square can contain a letter that does not appear in the rest of the array or is filled by one of
the $i-2s-r$ letters that go to themselves in the bottom $L$.
\end{enumerate}

\begin{center}
\includegraphics[width=6.0cm, height=3.5cm]{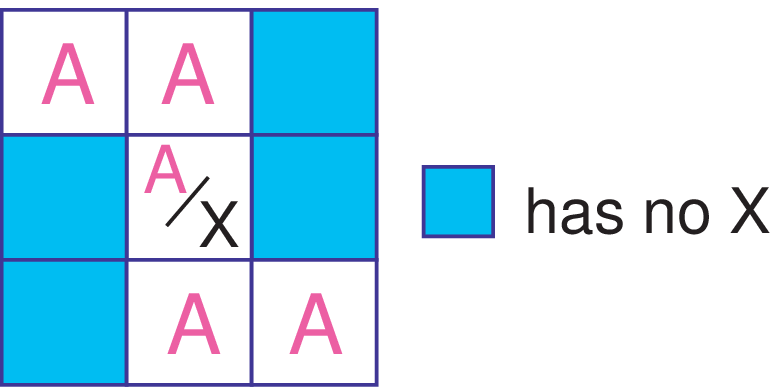}
\end{center}
\begin{center}
\small{Figure 4.2: The two ways to fill the central square.  The central square contains either new letter or a letter from in the upper $L$ that rotates to itself in the bottom $L$.}
\end{center}
\begin{theorem}
Let $m$ and $n$ be odd integers.  Let  $q = n[\frac{m}{2}] + [\frac{n}{2}]$.  The coefficient of $t^q$ in the expansion of $\exp(2e^t - 2 + t +\frac{1}{2}(e^t - 1)^2)$ is $R_{m,n}$. 
\end{theorem}
\newpage
\noindent {\bf Proof of Theorem 4.1:}\\
Let
\[R(t) = \sum_{\scriptstyle q,i,s,r = 0}^{\infty}\frac{\left\{\begin{array}{c} q\\i \end{array}\right\}t^qi!(i-2s-r+1)}{q!2^ss!r!(i-2s-r)!}\]
By using the techniques of Theorem 3.1, we obtain the desired result.\qquad $\Box$
\begin{remark}
We should note that the theorems in Sections 3 and 4 could be proven using the context of involutions acting on a finite set.  In particular, suppose the involution is acting on a finite set in manner
which provides $t$ matching pairs and $u$ fixed points.  Then, by the argument used to derive $H_{2m+1,n}$, the number of partitions of the set that are inequivalent under the involution is
$\exp(e^{y+x}+\frac{1}{2}e^{2y}-\frac{3}{2})$, where the powers of $x$ correspond to $u$ and the powers of $y$ correspond to $t$.  Thus, $H_{2m,n}$ is the case when $t = mn$ and $u = 0$, $H_{2m+1,n}$ is
the case when $t = mn$ and $u = n$, and $R_{m,n}$ (with $m, n$ odd) is the case when $t = q = n[\frac{m}{2}] + [\frac{n}{2}]$ and $u = 1$. 
\end{remark}

\section{Generating Function for Fully Symmetrical Letter Representation}
\noindent We define an $m\times n$ letter representation to be {\bf fully symmetric} if and only if it is fixed via horizontal \underline{and} vertical reflection.
 Let $S_{m,n}$ count the number of fully symmetrical letter
representations.  Although $S_{m,n}$ only appears in the proof of Theorem 2.1, we can apply the techniques of Sections 3 and 4 to determine the generating functions for $S_{m,n}$.  Since the
fully symmetrical letter representations are invariant under two symmetry transformations, the underlying $m\times n$ array is divided into four quadrants.  The upper left quadrant is arbitrarily
filled and the remaining three quadrants are completed by applying the two symmetry transformations to the upper left quadrant.  There are two types of letters in the upper left quadrant.  The first type
of letter is called a  singleton letter.  A {\bf singleton  letter} is a letter whose image, in the remaining three quadrants, is never another letter that appears in the upper left
quadrant.  There are five possible ways a letter can be a singleton.  These five ways are illustrated in Figure 5.1
\newpage
\begin{center}
\includegraphics[width=8.0cm, height=4.0cm]{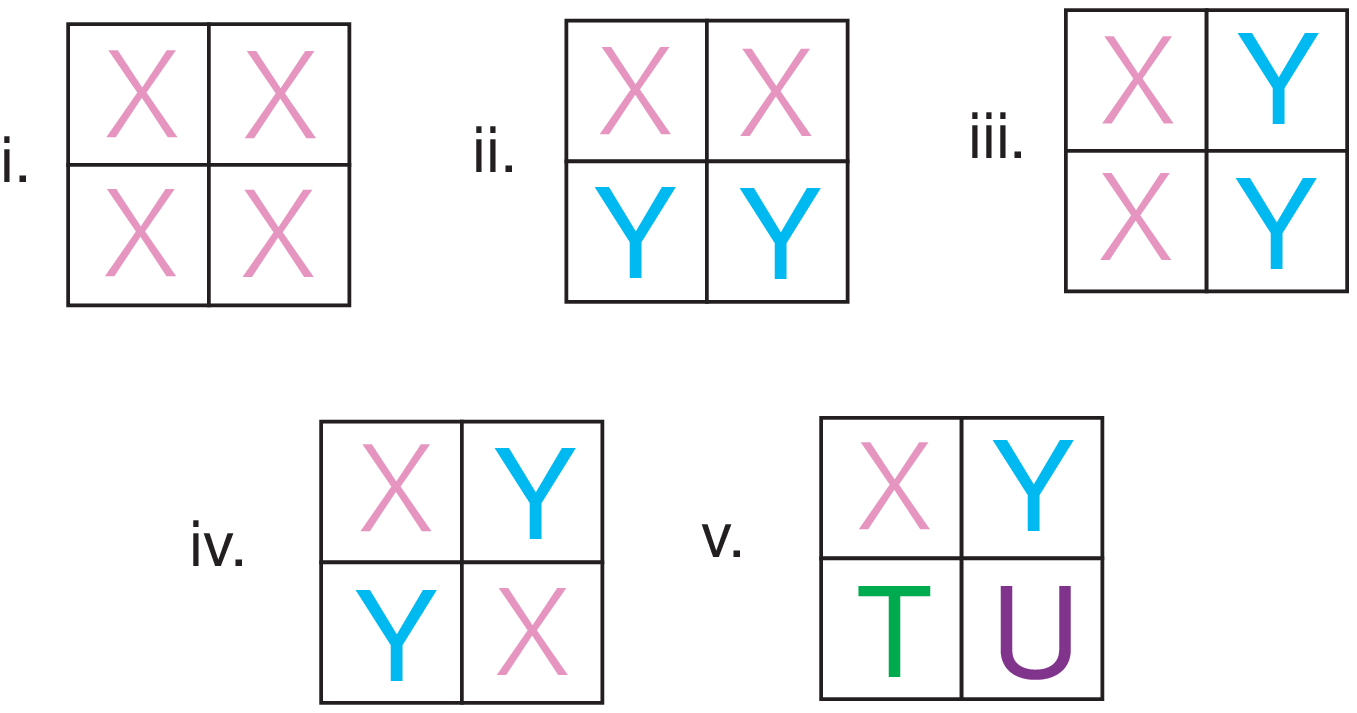}
\end{center}
\begin{center}
\small{Figure 5.1: The five ways a singleton letter can be transformed under symmetry.}
\end{center}
\noindent  The second type of letter present in the upper left quadrant can be considered to be part of a double pair.  A letter is part of {\bf double pair} when its image in one of the other
three quadrants is another letter orginally present in the upper left quadrant.  Figure 5.2 illustrates the six ways double pairs transform in a fully symmetrical manner.
\begin{center}
\includegraphics[width=6.0cm, height=10.0cm]{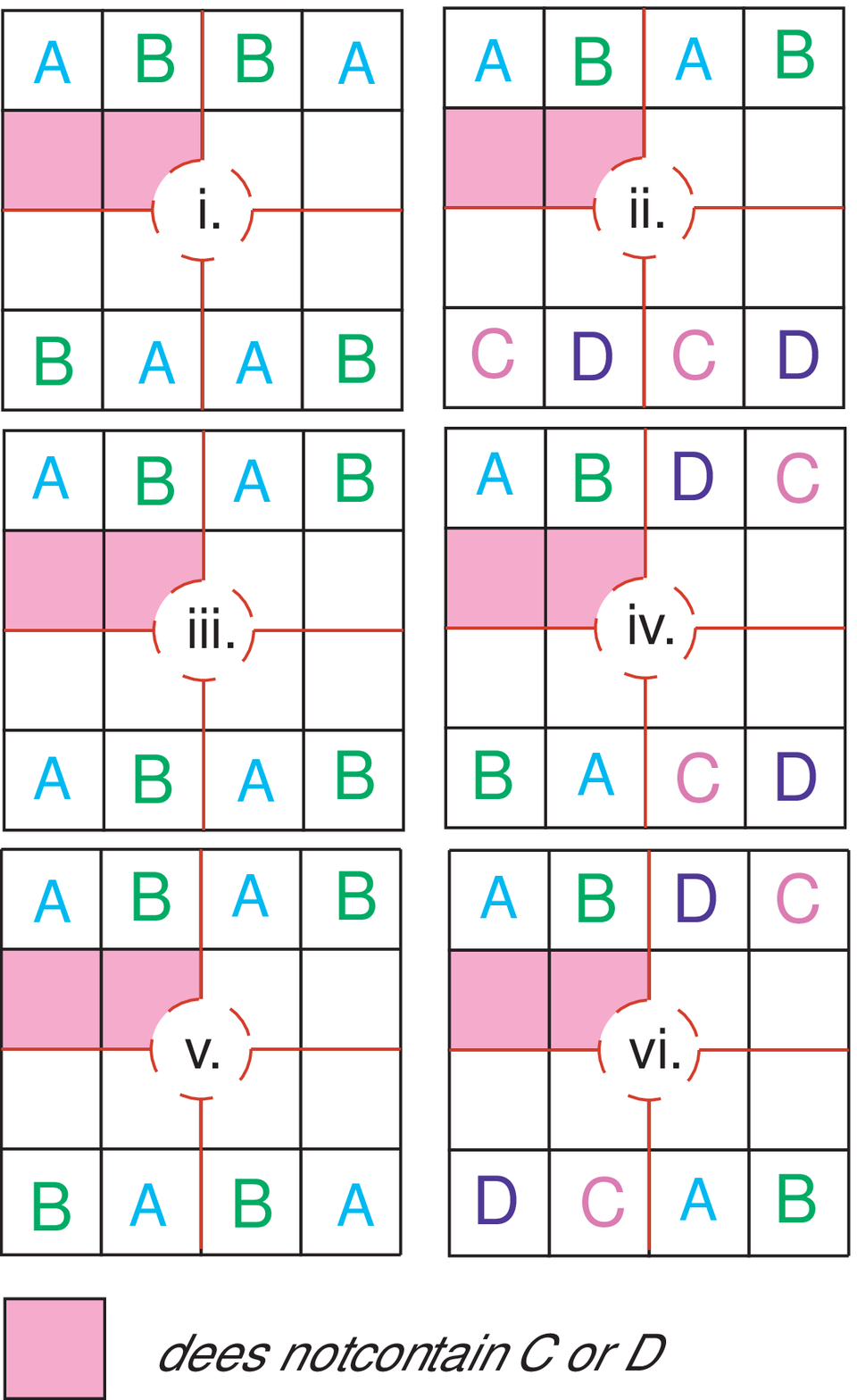}
\end{center}
\begin{center}
\small{Figure 5.2:  The six ways two letters, each appearing in the upper left quadrant, form double pairs. }
\end{center}
\noindent We begin by assuming the row and column dimensions are even.  The geometric construction depends on locating the singleton letters and the
double pairs present in the upper left quadrant.  In particular,   
\newpage
\[S_{2m,2n} = \sum_{j=1}^{mn}\sum_{s=0}^{[\frac{j}{2}]}\frac{\left\{\begin{array}{c} mn\\j \end{array}\right\}6^s5^{j-2s}j!}{2^s s!(j-2s)!}\]
\noindent where
\begin{enumerate}
\setlength{\itemsep}{-\itemsep}
\item[1.] $j$ counts the letters in the upper left quadrant. 
\item[2.] $s$ counts the number of these $j$ letters that are interchanged under symmetry operations.
\end{enumerate}
\noindent The proof of Theorem 5.1 is similiar to that of Theorem 3.1 and will be omitted.
\begin{theorem}
The coefficient of $t^{mn}$ in the expansion of $\exp(5(e^t - 1)e^{3(e^t + 1)^2})$ is $S_{2m,2n}$.  
\end{theorem}
\noindent Next we calculate $S_{2m,2n+1}$.  By symmetry, the calculation for $S_{2m,2n+1}$ determines $S_{2n+1,2m}$.  In this particular case, we decompose the $2m\times (2n+1)$ array into four quadrants
and a middle column and fill the upper left quadrant with an arbitrary letter configuration.  Using symmetry, we complete the remaining three quadrants as illustrated by Figures 5.1 and 5.2.  However, we
need to determine whether a singleton letter or letter from a double pair may occur in the middle column.  In particular, the two singleton types illustrated by Diagrams i and ii in Figure 5.1 and the
double pair illustrated by Diagram i of Figure 5.2 may occur in the middle column.  The middle column is completed as follows.  The first $m$ squares of the middle column, called the {\bf top half} of the
middle column, are filled by letters that appeared in the quadrants or by new letters.  Then, using horizontal reflection, the corresponding positions in the lower half of the middle column are filled. 
Hence
\[S_{2m,2n+1} = \sum_{i=1}^{mn}\sum_{F=0}^{[\frac{i}{2}]}\sum_{M=0}^{[\frac{i-2F}{2}]}\,\sum_{f=0}^{i-2F-2M}\,\sum_{p=0}^{i-2F-2M-f}\sum_{l=0}^{m}\,\sum_{k=min(1,m-l)}^{m-l}\sum_{t=0}^{[\frac{k}{2}]}\]
\[
\frac{\left\{\begin{array}{c} mn\\i \end{array}\right\}i!\,5^F3^f2^{-F-M+k-3t}m!(2i-2f-p-4F-2M)^l\left\{\begin{array}{c} m-l\\k\end{array}\right\}k!}{F!M!f!p!(i-2F-2M-f-p)!\,l!\,(m-l)!(k-2t)!\,t!}\]
where
\begin{enumerate}
\setlength{\itemsep}{-\itemsep}
\item[1.] $i$ counts the letters in the upper left quadrant.  
\item[2.] $F$ counts the double pairs of the upper left hand quadrant that cannot occur in the middle column.
\item[3.] $M$ counts the double pairs of the upper left quadrant that can occur in the middle column.
\item[4.] $f$ counts the singleton letters of the upper left quadrant that cannot occur in the middle column.
\item[5.] $p$ counts the letters in the upper left quadrant that are the singleton letters illustrated by Diagram i of Figure 5.1
\item[6.]$l$ counts the squares in the top half of the middle column that are filled by letter that appeared in the quadrants
\item[7.] $k$ counts the letters in the top half of the middle column that do not appear in the quadrants.
\item[8.] $t$ counts the interchanges that occur as the $k$ letters are reflected to the bottom half of the middle column.
\end{enumerate}
\begin{theorem}
The coefficient of $y^mx^{mn}$ in the expansion of $\exp(\frac{1}{2}[e^x(-4 + 2e^y) + e^{2x}(e^{2y} + 5) - 4])$ is $S_{2m,2n+1}$.
\end{theorem}
\noindent To prove Theorem 5.2, we make use of the following version of Taylor's Theorem.
\begin{remark}
\[e^{a\frac{d^n}{dw^n}}e^{bw} = e^{ab^n}e^{bw}\]
\end{remark}
\noindent {\bf Proof of Theorem 5.2}
Let $q = mn$. Then,
{\small
\begin{align*}
S(x,y) &= \sum_{\scriptstyle i,F,M,f,p, \atop \scriptstyle m,l,k,t,q = 0}^{\infty}
\frac{\left\{\begin{array}{c} q\\i \end{array}\right\}\left\{\begin{array}{c}
m-l\\k\end{array}\right\}i!\,k!(2i-2f-p-4F-2M)^l5^F3^f2^{-F-M+k-3t}x^qy^m}{F!M!f!p!(i-2F-2M-f-p)!\,l!\,(m-l)!\,t!\,(k-2t)!q!}\\
&= \sum_{\scriptstyle i, F, M, f, \atop \scriptstyle p,l,k,t = 0}^{\infty}
\frac{(e^x-1)^i(e^y-1)^k5^F3^f2^{-F-M+k-3t}(2i-2f-p-4F-2M)^ly^l}{F!M!f!p!(i-2F-2M-f-p)!\,l!\,(k-2t)!t!}\\
&= \sum_{\scriptstyle i, F, M, f,\atop \scriptstyle p,l,k,t = 0}^{\infty}\frac{(e^x-1)^i (e^y-1)^k 5^F3^f2^{-F-M-t}(2i-2f-p-4F-2M)^ly^l[dz]^{2t}z^k}{F!M!f!p!(i-2F-2M-f-p)!\,l!\,k!\,t!}\Bigg|_{z=2}\\
&= \sum_{\scriptstyle i, F, M, f,\atop \scriptstyle p,l,k =0}^{\infty}\frac{(e^x-1)^i(e^y-1)^k5^F3^f2^{-F-M}(2i-2f-p-4F-2M)^ly^le^{\frac{(dz)^2}{2}}z^k}{F!M!f!p!(i-2F-2M-f-p)!\,l!\,k!}\Bigg|_{z=2}
\end{align*}}
\noindent Due to Remark 5.1, we obtain
\begin{align*}
S(x,y)&= \sum_{\scriptstyle i, F, M,f,p,l = 0}^{\infty}\frac{(e^x-1)^i 5^F3^f2^{-F-M}(2i-2f-p-4F-2M)^ly^l e^{\frac{1}{2}(e^{2y}+2e^y-3)}}{F!M!f!p!(i-2F-2M-f-p)!\,l!}\\
&= \sum_{\scriptstyle i, F, M, f,p = 0}^{\infty}\frac{(e^x-1)^i 5^F3^f2^{-F-M}e^{y(2i-2f-p-4F-2M)}e^{\frac{1}{2}(e^{2y}+2e^y-3)}}{F!M!f!p!(i-2F-2M-f-p)!}
\end{align*}
\noindent  Let $M = u - F$.  The above line becomes
\begin{align*}
S(x,y)&= \sum_{\scriptstyle i, F, u,f,p = 0}^{\infty}\frac{[(e^x-1)e^{2y}]^i[3e^{-2y}]^f[5e^{-2y}]^F[e^{-y}]^p[2e^{2y}]^{-u}e^{\frac{1}{2}(e^{2y}+2e^y-3)}} {F!(u-F)!f!p!(i-2u-f-p)!}\\
&= \sum_{\scriptstyle i, F, u = 0}^{\infty}\frac{[(e^x-1)e^{2y}]^i[5e^{-2y}]^F[1+3e^{-2y}+e^{-y}]^{i-2u}[2e^{2y}]^{-u}e^{\frac{1}{2}(e^{2y}+2e^y-3)}}{F!(u-F)!(i-2u)!}\\
&= \sum_{\scriptstyle i, u =0}^{\infty}\frac{[(e^x-1)e^{2y}]^i[1+5e^{-2y}]^u[1+3e^{-2y}+e^{-y}]^{i-2u}[2e^{2y}]^{-u}e^{\frac{1}{2}(e^{2y}+2e^y-3)}} {u!(i-2u)!}\\
&= \exp(\frac{1}{2}[e^x(-4+2e^y)+(e^{2y}+5)e^{2x}-4])\qquad \Box
\end{align*}

\noindent The final sum of this section counts $S_{(2m+1)\times (2n+1)}$.  In the case, the geometric decomposition involves four quadrants, a central
square, and a middle cross, where a {\bf middle cross} is the union of the middle column and middle row minus the central square.   Thus, a middle cross consists of the top and bottom halves of the middle
column along with the left and right halves of the middle row.  In this situation, {\bf the top half of the middle column} is the $m$ squares above the central square, and {\bf the left half of the middle
row} is the $n$ squares to the left of the central square.\\\\ 
\noindent For this situation, we will not begin by analyzing the components in the upper left quadrant.  Instead, we will
work first with the middle cross and then extend outward to the quadrants.  In particular, we begin by filling in the top half of the middle column in an
arbitrary manner.  Next, using the horizontal reflective symmetry, we determine bottom half of the middle column.  A letter from the top half, when reflected to the
bottom, can either go to itself, go to an entirely new letter, or become interchanged with another letter present in the top half.  In all three cases, these top half letters (and possibly new bottom
half letters) may occur in the quadrants.  This situation is illustrated by Figure 5.3.
\begin{center}
\includegraphics[width=8.0cm, height=10.0cm]{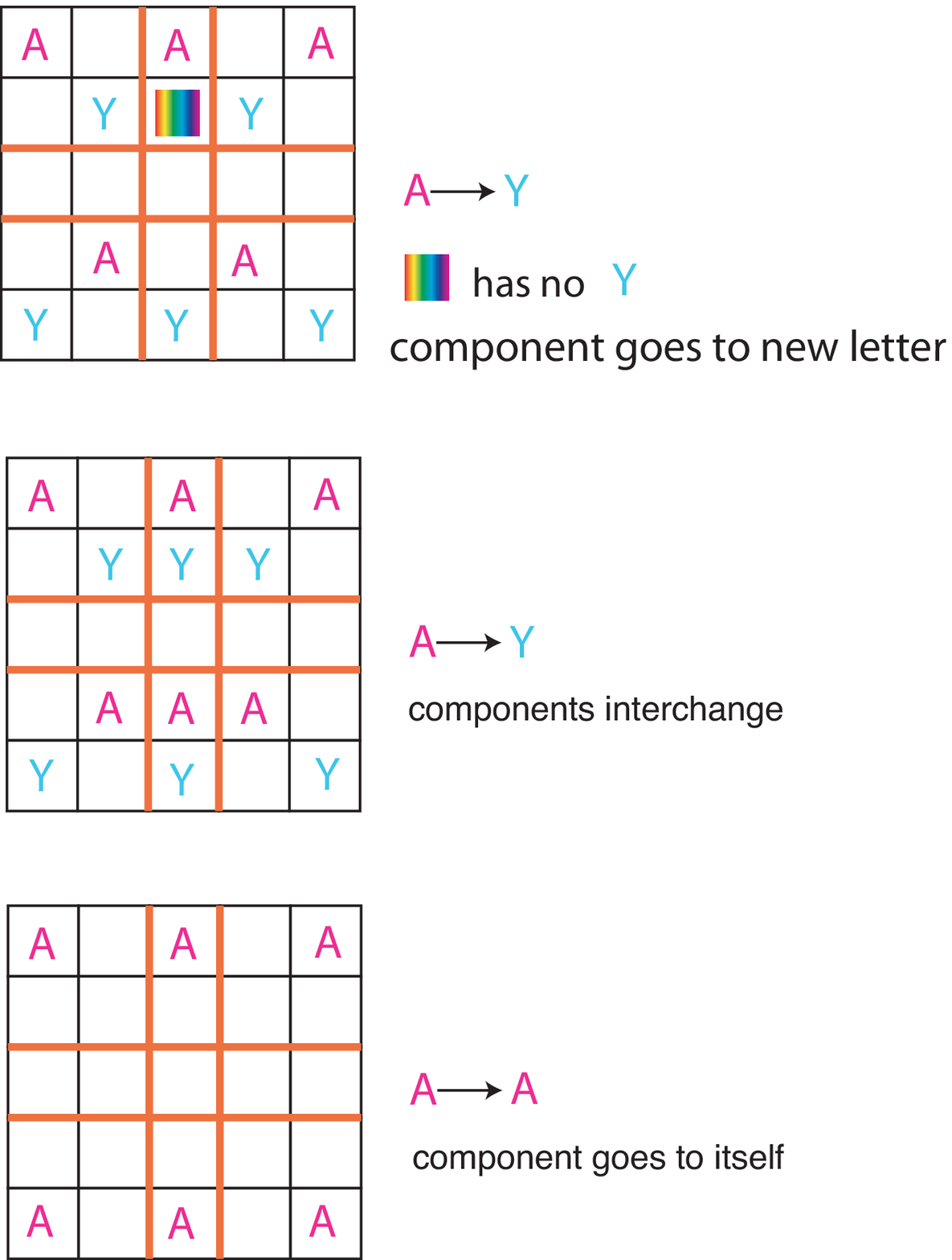}
\end{center}
\begin{center} 
\small{Figure 5.3: The three ways $A$, a letter from the top half of the middle column, is transformed, via horizontal reflection, to the bottom half of the middle column.  In all three cases, note that $A$
and its image has a unique way of appearing in the quadrants.}
\end{center}
\noindent The next step is to determine the middle row.  This is done by first filling the left half of the middle row with
letters occuring in the middle column.  The only such middle column letters that may be used are those letters which are fixed under horizontal reflection.  To fill in the remaining spaces of the middle
row minus the central square, we repeat the procedure used to determine the middle column.  The only difference is that we use a vertical reflection.  Once again, these middle row letters may appear in the
quadrants.  To visualize this situation, simply rotate the diagrams in Figure 5.5 by 90 degrees.  Next, we are in a position to determine the quadrants.  This is done by first placing the cross letters
throughout the quadrants.  Then, for the remaining spaces of the quadrants, we repeat the geometric analysis of singletons and double pairs used for determining $S_{2m,2n}$.  The last step is to fill the
central square with an entirely new letter, a singleton from the quadrant that is fixed under symmetry, or a letter from the middle cross that is also fixed under symmetry.\\\\
\noindent Following the steps outlined above, we arrive at a summation of 22 variables.  The derivation of the exponential generating function utilizes the techniques of Theorem 5.2.  Details are
available, upon request, from the author.
\begin{theorem}
The coeffiecient of $y^mx^nl^{mn}$ in the expansion of\\ $\exp(l + x + y)\exp(2e^{2l} + \frac{1}{2}e^{2l+2x}+\frac{1}{2}e^{2l+2y}+e^{l+x+y}-2e^l-2)$ is $S_{2m+1,2n+1}$.
\end{theorem}
\section{Further Research}
\noindent By using a particular decomposition of the Bell Numbers [4],[7],[8] and applying various symmetry transformations to $m\times n$ rectangular arrays, we are able to calculate
the number of $m\times n$ letter representations modulo $D_2$ symmetry.  The number of $m\times n$ letter representations only provides a lower bound for the basis size of the transition
matrix.  In order to find an upper bound for the cardinality of the basis, it is necessary to extrapolate the results of this paper to the case of $m\times n$ arrays of circled letters (i.e. word
representations of Section 1).  Currently, we are formulating the generating functions used in the enumeration of $m\times n$ arrays of circled letters.\\\\  
\noindent Another promising avenue of research involves square arrays, where $m = n$.  Since a square array has $D_4$ symmetry, the formula to enumerate $m\times m$ letter representations modulo symmetry
will be more complicated than the formula provided by Theorem 2.1.  We are in the process of calculating the
generating functions necessary in the enumeration of $m\times m$ letter representations.  The ultimate goal will be to adapt the geometric constructions
used in the enumeration of $m\times m$ letter representations to the case of $m\times m$ word representations.  Once this is accomplished, an upper bound for the basis size
of the associated transition matrix will be completely determined for all integers $m$ and $n$.\\\\
\noindent Finally, a different type of research would involve exploring the connections between $m\times n\times p$ proper arrays and percolation theory.  At the present time, we have not explored
the connection in any depth but realize that the stochastic and probabilisitic techniques of percolation theory could, when applied to the representation of an $m\times n\times p$ proper array as a bond
percolation on $Z^3$ with an open cluster at the origin (see Section 0.3), give rise to a whole new category of results.}
\section*{Acknowledgments}
\noindent The author would like to thank Dr. Bruce Sagan and Dr. Robert Sulanke for their suggestions regarding the Introduction.  The author also thanks Dr. Harris Kwong for his time and help in editing
the previous drafts of this paper.
\newpage
\section*{Appendix A: Numerical Data}
\noindent The following table provides, for small integer values of $m$ and $n$, numerical values for $H_{m,n}, V_{m,n}, R_{m,n}$, and $S_{m,n}$.  All
the values came from the generating functions given by Theorems 3.1, 3.2, 4.1, 5.1, 5.2, and 5.3 and were verified by a short Maple program the author created.
\begin{center}
\begin{tabular}{||l|c|c|c|c||}
\hline
$m\times n$ & $H_{m,n}$ & $V_{m,n}$ & $R_{m,n}$ & $S_{m,n}$\\ \hline
$2\times 3$ & 31 & 12 & 31 & 13\\ \hline
$2\times 4$ & 164 & 164 & 164 & 36\\ \hline
$2\times 5$ & 999 & 999 & 999 & 107 \\ \hline
$3\times 2$ & 12 & 31 & 31 & 13 \\ \hline
$3\times 4$ & 339 & 6841 & 6841 & 469\\ \hline
$3\times 5$ & 2210 & 51790 & 127643 & 3835 \\ \hline
\end{tabular}
\end{center}
\begin{center}
\small{Table 2:  Numerical Data for certain $m\times n$ letter representations}
\end{center}
\newpage
\section*{References}
$[1]$\hspace{0.2cm} J. Quaintance ``$m\times n$ Proper Arrays: Geometric and Algebraic Methods of Classification'', Ph.D. Dissertation, University of Pittsburgh, August 2002\\\\
$[2]$\hspace{0.2cm} J. Quaintance ``$n\times m$ Proper Arrays: Geometric Constructions and the Associated Linear Cellular Automata'', {\it Maple Summer Workshop 2004 Proceedings}\\\\
$[3]$\hspace{0.2cm} J. Quaintance ``Combinatoric Enumeration of the Geometric Classes Associated with $n\times p$ Proper Arrays'', preprint January 2004\\\\ 
$[4]$\hspace{0.2cm} N. J. A. Sloane, ``The On-Line Encyclopedia of Integer Sequences'',\\ http://www.research.att.com/~njas/sequences/ \\\\
$[5]$\hspace{0.2cm} S. Lo and M. Monagan, `` A Modular Algorithm for Computing the Characteristic Polynomial of an Integer Matrix in Maple '', {\it Maple Summer Workshop 2005 Proceedings}\\\\
$[6]$ \hspace{0.2cm} R. Simion, `` Noncrossing partitions '', {\it Discrete Mathematics}, {\bf 217} (2000), 367-409\\\\
$[7]$\hspace{0.2cm} K. Yoshinaga and M. Mori, `` Note on an Exponential Generating Function of Bell Numbers '', {\it Bull. Kyushu Inst. Tech.}, {\bf 24} (1976), 23-27\\\\
$[8]$\hspace{0.2cm} D. Branson, ``Stirling Numbers and Bell Numbers: Their Role in Combinatorics and Probability '', {\it Math. Scientist}, {\bf 25} (2000), 1-31\\\\
$[9]$\hspace{0.2cm} B. Cipra, `` An Introduction to the Ising Model '', {\it The American Mathematical Monthly}, {\bf 94} (1987), No. 10, 937-959\\\\
$[10]$\hspace{0.2cm} G. Grimmett, {\it Pecolation} Second Edition, Springer-Verlag, 1991\\\\
$[11]$\hspace{0.2cm} D. Stauffer and A. Aharony, {\it Introduction to Percolation Theory} Second Edition, Taylor and Francis, 1992\\\\
$[12]$\hspace{0.2cm} R. Stanley, {\it Enumerative Combinatorics, Volume 1} New York: Wadsworth and Brook/Cole, 1986 \\\\
\newpage
\end{document}